\documentclass[12pt]{article}
\usepackage{amssymb}

\def\header{A.Miller\hfill Long Borel Hierarchies \hfill} 
\def\al{\alpha}
\def\be{\beta}
\def\de{\delta}
\def\ga{\gamma}
\def\ka{\kappa}
\def\om{\omega}
\def\si{\sigma}

\def\la{\langle}
\def\ra{\rangle}

\def\bpi(#1){{\bf \Pi}^0_{#1}}
\def\bsi(#1){{\bf \Sigma}^0_{#1}}

\def\proof{\par\noindent Proof\par\noindent}
\def\qed{\par\noindent QED\par}

\def\aa{{\mathcal A}}
\def\bb{{\mathcal B}}
\def\cc{{\mathcal C}}
\def\ff{{\mathcal F}}
\def\gg{{\mathcal G}}
\def\mm{{\mathcal M}}
\def\nn{{\mathcal N}}
\def\qq{{\mathcal Q}}

\def\pow{{\mathcal P}}
\def\fix{{\rm fix}}
\def\col{{\mathbb Coll}}
\def\comp(#1){\sim #1}
\def\rmor{\mbox{ or }}
\def\clash{\#}
\def\sde{{s\hat{\phantom{a}}\la\de\ra}}
\def\tde{{t\hat{\phantom{a}}\la\de\ra}}
\def\sin{{s_i\hat{\phantom{a}}\la n\ra}}
\def\scatn{{s\hat{\phantom{a}}\la n\ra}}
\def\alcats{{\la\al\ra\hat{\phantom{a}}s}}
\def\Ga{\Gamma}
\def\vleq{\unlhd}
\def\child{{\rm Child}}
\def\name#1{\stackrel{\circ}{#1}}
\def\cl{{\rm cl}}
\def\rank{{\rm rank}}
\def\cof{{\rm cof}}
\def\poset{{\mathbb P}}
\def\leaf{{\rm Leaf}}
\def\dom{{\rm dom}}
\def\rmand{\mbox{ and }}
\def\rmiff{\mbox{ iff }}
\def\group{{\mathcal H}}
\def\su{\subseteq}
\def\supp{{\rm supp}}
\def\sm{{\setminus}}
\def\st{\;:\;} 
\def\forces{{\;\Vdash}}

\def\res{\upharpoonright}
\def\pr{\prime}
\def\rat{{\mathbb Q}}

\newtheorem{theorem}{Theorem} 
\newtheorem{lemma}[theorem]{Lemma}
\newtheorem{cor}[theorem]{Corollary}
\newtheorem{question}[theorem]{Question}
\newtheorem{conj}[theorem]{Conjecture}
\newtheorem{prop}[theorem]{Proposition}

\begin{document}

\begin{center}
{\large Long Borel Hierarchies}
\end{center}

\begin{flushright}
Arnold W. Miller
\footnote{
Thanks to the University of Florida Mathematics Department for their support and
especially Jindrich Zapletal,  William Mitchell, Jean A. Larson, and Douglas
Cenzer for inviting me to the special year in Logic 2006-07 during which this
work was done.  Also I would like to thank P\'eter Komj\'ath for telling me
about this problem. \par Mathematics Subject Classification 2000: 03E15; 03E25;
03E35 \par Keywords: Axiom of Choice, Borel Hierarchies, Meager Ideal, Countable
Unions. }
\end{flushright}

\def\address{\begin{flushleft}
Arnold W. Miller \\
miller@math.wisc.edu \\
http://www.math.wisc.edu/$\sim$miller\\
University of Wisconsin-Madison \\
Department of Mathematics, Van Vleck Hall \\
480 Lincoln Drive \\
Madison, Wisconsin 53706-1388 \\
\end{flushleft}}

\begin{center}  Abstract  \end{center}
\begin{quote}
We show that it is relatively consistent with ZF that the
Borel hierarchy on the reals has length $\om_2$.
This implies that $\om_1$ has countable cofinality, so
the axiom of choice fails very badly in our model.
A similar argument produces models of ZF in which the Borel
hierarchy has length any given limit ordinal less than $\om_2$, 
e.g., $\om$ or $\om_1+\om_1$.  
\end{quote}

\begin{center}
Introduction
\end{center}

In this paper we do not assume the axiom of choice, not even in the
form of choice functions for countable families.
Define the classical Borel families, $\bpi(\al)$ and $\bsi(\al)$, of subsets
of $2^\om$ for any ordinal $\al$ as usual:

\begin{enumerate}
\item  $\bsi(0)=\bpi(0)=$clopen subsets of $2^\om$,
\item  $\bpi(<\al)=\cup_{\be<\al}\bpi(\be)$,
$\;\;\;\;\;\;\;\;\;\;$  $\bsi(<\al)=\cup_{\be<\al}\bsi(\be)$,
\item $\bsi(\al)=\{\cup_{n<\om}A_n\st (A_n:n<\om)\in
(\bpi(<\al))^\om\}$, and
\item $\bpi(\al)=\{\cap_{n<\om}A_n\st (A_n:n<\om)\in
(\bsi(<\al))^\om\}$.
\end{enumerate}

\noindent It follows immediately from these definitions that for all $\al<\be$
$$\bpi(\al)\cup\bsi(\al)\su \bpi(\be)\cap\bsi(\be).$$
Hence for limit ordinal $\al$ we have that $\bpi(<\al)=\bsi(<\al)$.
It is also true by DeMorgan's Laws that
$$\bpi(\al)=\{(2^\om\sm X)\st X\in\bsi(\al)\}.$$

The family of Borel subsets of $2^\om$ is the smallest family of sets
containing the clopen sets and closed under countable unions
and countable intersections.
Equivalently,
Borel$=\bsi(<\infty)=\bpi(<\infty)$ where
$$\bsi(<\infty)=\cup\{\bsi(\al)\st \al \mbox{ an ordinal }\}
\rmand
\bpi(<\infty)=\cup\{\bpi(\al)\st \al \mbox{ an ordinal }\}.$$

Let us call the least $\al$ such that Borel $=\bsi(<\al)$
the length of the Borel hierarchy. The length cannot be
$\infty$ since then
there would be a map from the power set of $2^\om$ onto the class of all
ordinals.

It is a classical Theorem of Lebesgue 1905 \cite{leb} (see Kechris
\cite{kechris}) that assuming the axiom of choice for countable families, the
length of the Borel hierarchy is $\om_1$. To see that it has height at least
$\om_1$, he shows that $\bsi(\al)\neq\bpi(\al)$  for all $\al$ with
$1\leq\al<\om_1$. In the absence  of the axiom of choice this may fail.
Feferman and Levy 1963 (see Jech \cite{jech}) showed that it is relatively
consistent with ZF that $2^\om$ is the countable union of countable sets.  This
implies that every subset of $2^\om$ is a countable union of countable sets.
Hence in the Feferman-Levy model every subset of  $2^\om$ is Borel and the
Borel hierarchy has finite length. In their model the $\bsi(2)$ sets are  not
closed under countable unions.

The place in the Lebesgue proof which goes wrong is in the construction
of a universal set for each Borel class.  This requires choosing codes
for Borel sets.

Since there is a map of $2^\om$ onto $\om_1$, it also true
in the Feferman-Levy
model that $\om_1$ has cofinality $\om$. In fact, in their model
$\om_1=\aleph_{\om+1}^L$.  Lebesgue also needed the axiom of
choice to see that $\om_1$ is a regular cardinal, and therefor each Borel set
will appear at a countable level of the Borel hierarchy, i.e.
$\bsi(<\om_1)=\bpi(<\om_1)=$ Borel.

P\'eter Komj\'ath asked if it is possible for the Borel hierarchy to have length
greater than $\om_1$ in some model of ZF.  We show that it can be.
This is the main result of our paper.

\begin{quote}
{\bf Theorem \ref{thm1}. }
{\it It is relatively consistent with ZF that the Borel hierarchy on
$2^\om$ has length $\om_2$, i.e., the least $\al$ such that
$\bsi(<\al)$ is the
family of all Borel sets is $\al=\om_2$.}
\end{quote}

Our model will be a symmetric submodel $\nn$ of a generic
extension of the Feferman-Levy model $V$. 
In an inner model of our main model we will find models of 
ZF in which the Borel
hierarchy has length any given limit ordinal less than $\om_2$
(Theorem \ref{thm2}).
Using a model of Gitik \cite{gitik} in which
every cardinal is singular, we show 
that the Borel hierarchy can be arbitrarily high 
(Theorem \ref{gitikmodel}).

\begin{center}
Proof of Theorem \ref{thm1}.
\end{center}

The Feferman-Levy Model $V$ is described in Jech \cite{jech}.  The
ground model satisfies $V=L$, let us call it $L$.
In $L$ let $\col$ be the following version of the
Levy collapse of $\aleph_\om$:
$$\col=\{p:F\to \aleph_{\om}\st F\in[\om\times\om]^{<\om}\rmand
\forall (n,m)\in F\;\; p(n,m)\in\aleph_n\}$$
ordered by inclusion.

For any $n<\om$ let
$\col_n=\{p\st \dom(p)\subseteq n\times\om\}$
and for $G^c$ $\col$-generic over $L$ let $G^c_n=G^c\cap\col_n$.

The properties we will use of $V$ are summarized in the
next Lemma.

\begin{lemma}\label{F-L}
$L\su V\su L[G^c]$ and $G^c_n\in V$ for each $n$.
In $V$
$\pow(\om)$ is the countable union of countable sets, in fact,
$$ \pow(\om)\cap V=\bigcup_{n<\om}(L[G^c_n]\cap \pow(\om)).$$
More generally,
if $X\su Y\in L$ and $X\in V$, then for some $n<\om$ we have that
$X\in L[G_n]$.  It follows that
$\om_1^V=\aleph_\om^L$ and
$\om_2^V=\aleph_{\om+1}^L$ and is regular in $V$.
\end{lemma}

Working in $L$ construct a well-founded tree $T\su (\aleph_{\om+1})^{<\om}$.

\bigskip
\noindent  First we define:
\begin{enumerate}
\item For $s\in (\aleph_{\om+1})^{<\om}$ and $\de<\aleph_{\om+1}$,
$\;\;\sde$ is the finite sequence of length
$|s|+1$ which begins with $s$ and has one more
element $\de$.
\item For $s\in T$  $$\child(s)=\{\de\st \sde\in T\}.$$
\item For $s\in T$ 
$$\rank(s)=\sup\{\rank(\sde)+1\st \de\in\child(s)\}$$
Note that $\rank(s)=0$ for terminal nodes or leaves, $s\in\leaf(T)$.
\end{enumerate}

\bigskip
\noindent Then $T$ should have the following properties:
\begin{enumerate}
\item $\child(\la\ra)=\aleph_{\om+1}$ and
$\rank(\la\al\ra)=\al$
for each $\al<\aleph_{\om+1}$.
\item If $\rank(s)=\al+1$ a successor ordinal, then
$\{\de\st \sde\in T\}=\om$ and
$\rank(\scatn)=\al$ for all $n<\om$.
\item If $\rank(s)=\lambda$ a limit ordinal and $\cof(\lambda)=\om_n$, then
$\child(s)=\om_n$ and $\rank(\sde)$ for $\de<\om_n$ is
strictly increasing and (necessarily) cofinal in $\lambda$.
\end{enumerate}

It is easy to inductively construct such a $T$ in $L$.  Note
that in $V$ each $\om_n^L$ is countable, so except for the
root node $\la\ra$, $T$ is countably branching, i.e.,
$\child(s)$ is countable for every $s\in T$ except the root node.

\begin{enumerate}
\item For $\leaf(T)$ the terminal nodes of $T$,
define $\poset$ to be the set of finite partial functions
$p:F\to 2^{<\om}$ for $F\in [\leaf(T)]^{<\om}$ ordered by
$p\leq q$ iff $\dom(p)\supseteq \dom(q)$ and $p(s)\supseteq q(s)$ for
every $s\in \dom(q)$.  This is forcing equivalent to Cohen real forcing,
$FIN(\aleph_{\om+1},2)$.
\item For $\pi$ a permutation, define the support of $\pi$,
$$\supp(\pi)=\{t\in\dom(\pi)\st \pi(t)\neq t\}.$$
\item Let $\group$ be the group of automorphisms of $\poset$
which are induced by finite support permutations of $\leaf(T)$.
That is, $\pi\in\group$ iff there exists a finite support
permutation $\hat{\pi}:\leaf(T)\to\leaf(T)$
such that $\pi:\poset\to\poset$ is defined by
$$\dom(\pi(p))=\hat{\pi}(\dom(p)) \rmand \pi(p)(s)=p(\hat{\pi}(s)).$$
\item For any $r\in T$ put $\leaf(r)=\{t\in\leaf(T)\st r\su t\}$.
Note that $\leaf(s)=\{s\}$ for $s\in\leaf(T)$.
\item For any $s\in T\sm\leaf(T)$ define
$$H_s=\{\pi\in\group\st\hat{\pi}(\leaf(\sde))=\leaf(\sde)
\mbox{ for all } \de\in\child(s)\}.$$
\item For any $t\in\leaf(T)$ define
$H_t=\{\pi\in\group\st \hat{\pi}(t)=t \}$.
\item Let $\ff$ be the filter of subgroups of $\group$ which are
generated by the $H_s$', i.e., $H\in\ff$ iff
there is a finite $Q\su T$ with
$$H_Q\su H\su \group\mbox{ where } H_Q=^{def}\cap\{H_s\st s\in Q\}.$$
\end{enumerate}

Note that we defined $H_t$ for $t\in\leaf(T)$ just for
convenience of notation,
since if $\scatn=t$, then $H_s\su H_t$.

\begin{lemma}
The filter of subgroups $\ff$ is normal, i.e., for any $\pi\in \group$
and $H\in\ff$, we have that $\pi^{-1} H\pi\in \ff$.
\end{lemma}

\proof
Fix $\pi\in\group$ and $Q\su T$ finite with $H_Q\su H$.
Let $R$ be a
finite superset of $Q$ which contains the support
of $\hat{\pi}$.  We claim that $\pi H_{R}\pi^{-1}= H_{R}$.
This follows from the fact that for any $\si\in H_{R}$ the
support of $\hat{\si}$ is disjoint from the support of
$\hat{\pi}$ and so $\pi\si\pi^{-1}=\si$.

It follows that
$$ \pi H_{R}\pi^{-1}= H_{R}\su H_Q\mbox{ implies }H_{R}\su \pi^{-1}H_Q\pi
\su \pi^{-1} H\pi$$
and hence $\pi^{-1} H\pi$ is in $\ff$.
\qed

Let $G$ be $\poset$-generic over $V$ and let $\nn$ be the
symmetric model determined by $\group$ and $\ff$.

\begin{lemma}\label{regular}
$\om_1^V=\om_1^{\nn}$, $\om_2^V=\om_2^{\nn}$,
and $\om_2^\nn$ remains regular in $\nn$.
\end{lemma}
\proof
It is enough to verify that this is true for
$V[G]$ in place of $\nn$, since
$$V\su\nn\su V[G]$$
This would seem obvious since $\poset$ is forcing equivalent to
the poset of the finite partial functions, $FIN(\ka,2)$,
where $\ka$ is $\om_2^V=\aleph_{\om+1}^L$.
If $V$ were a model of the axiom of choice, then we would know
that forcing with $\poset$ cannot collapse cardinals.

First we verify that $\om_1^V=\aleph_{\om}^L$ is
not collapsed in $V[G]$.
Working in $V$, suppose
for contradiction there exists $p_0\in \poset$ and
a name $\tau$ such that
$$p_0\forces \tau:\om\to\aleph_{\om}^L\mbox{ is onto}.$$
Define
$$A=\{(p,n,\be)\in \poset\times \om\times \aleph_{\om}^L
\st p\leq p_0\rmand p\forces
\tau(n)=\check{\be}\}$$
Note that for any $(p,n,\be),(q,n,\ga)\in A$ that if
$\be\neq\ga$, then $p$ and $q$ are incompatible.

The set $A$ is a subset of a set in $L$, so
it follows from 
Lemma \ref{F-L}
that there exist $k<\om$ such that
$A\in L[G^c_k]$.  In $L[G^c_k]$, $\om_1$ is $\aleph_{k+1}^L$.
Since $L[G^c_k]$ is a model of the axiom of choice,  the range
of $A$, i.e., $\{\al \st \exists p,n\;\;(p,n,\al)\in A\}$,
cannot even cover $\aleph_{k+1}^L$.

Now suppose in $V$
$$p_0\forces \tau:\om\to\aleph_{\om+1}^L \mbox{ is cofinal.}$$
Define $A$ similarly and suppose $A\in L[G^c_k]$.
Then since $\om_2^V=\aleph_{\om+1}^L=\aleph_{\om+1}^{L[G^c_k]}$
it follows that the range of $A$ cannot be cofinal
in $\om_2^V=\aleph_{\om+1}^L$.
This shows that the cofinality of $\om_2$ is $\om_2$ in $V[G]$
and hence it is not collapsed and it remains
regular.\footnote{An alternative proof for $\om_2$ regular
in $V$ is to note that it is $\om_1$ in the model $L[G^c]$. 
Since $L[G^c]$
is a model of ZFC forcing with $FIN(\ka,2)$ cannot
collapse $\om_1$.  The proof of Theorem \ref{gitikmodel}  
has an alternative argument for showing that 
cardinals are not collapsed in $\nn$.}

\qed

\begin{question}
Can there be a model of ZF in which for some $\ka$ forcing with
$FIN(\ka,2)$ collapses a cardinal?
\end{question}

For each $t\in\leaf(T)$ let $x_t\in 2^\om$ be the Cohen real attached to
$t$ which is
determined by $G$, i.e.,
$$x_t=\cup\{p(t)\st t\in\dom(p) \rmand p\in G\}.$$

For each
$s\in T$ define
$$A_s=\{x_t\st t\in\leaf(s)\}.$$
So $A_{\la\ra}$ is the set of all Cohen reals.
Working in $\nn$ for each ordinal $\al$ define the family
$\aa_\al$ inductively as follows:
\begin{enumerate}
\item $\aa_0$ is the set of finite subsets of $2^\om$, i.e.
$\aa_0=[2^\om]^{<\om}$,
\item $\aa_{<\al}=\cup_{\be<\al}\aa_\be$,
\item $\aa_\al=\{\cup_{n<\om} X_n\st (X_n:n<\om)\in (\aa_{<\al})^\om\}$
\end{enumerate}

\begin{lemma} \label{easy}
For each $s\in T$ the set $A_s$ is in $\nn$. For each $s\in T$ (except
the root node) $\;\;\;A_s\in\aa_\al$ where $\rank(s)=\al<\om_2$.
\end{lemma}
\proof
If $s\in \leaf(T)$, then the name of $x_s$:
$$\name{x}_s=\{(p,\check{\la n,i\ra})\st p\in\poset,\;\; p(s)=\si,
\rmand \si(n)=i\}$$
is fixed by all $\pi\in H_s$.
For any $s\in T$ the set
$A_s=\{x_t: t\in\leaf(s)\}$
has the name $\name{A_s}=\{(1,\name{x}_t)\st t\in\leaf(s)\}$
which is fixed by $H_s$.

Fix $s\in T$ with $\rank(s)=\al<\om_2^\nn$ and assume by induction
that for every $\de\in\child(s)$ that $A_{\sde}\in\aa_{<\al}$.  Then
$H_s$ fixes each $\name{A}_{\sde}$
for $\de\in\child(s)$ and so it fixes a name for the sequence
$\la A_{\sde}:\de\in\child(s)\ra$.
So this sequence is in $\nn$.  Since $\child(s)$ is countable in $V\su\nn$,
we see that $A_s\in\aa_\al$.
\qed

The elements of $\aa_\al$ are Borel sets, since finite sets are closed.
Similarly in the model $\nn$ define
\begin{enumerate}
\item $\mm_0$ to be the nowhere dense subsets of $2^\om$,
i.e., sets whose closure has no interior,
\item $\mm_{<\al}=\cup_{\be<\al}\mm_\be$
\item $\mm_\al=\{\cup_{n<\om} X_n\st (X_n:n<\om)\in (\mm_{<\al})^\om\}$
\end{enumerate}

Note that $\aa_\al\su\mm_\al$ since finite sets are nowhere dense.   The
following Lemma is proved by induction on $\al$ and is also true for $\aa_\al$.

\begin{lemma}
For any ordinal $\al$ the family $\mm_\al$ is closed under finite unions and
subsets, i.e., if $X,Y\in\mm_\al$, then $X\cup Y\in\mm_\al$ and
if $X\su Y\in\mm_\al$, then $X\in\mm_\al$.
\end{lemma}
\proof
Left to reader.
\qed

The usual clopen basis for $2^\om$ consists of sets of
the form
$$[\si]=\{x\in 2^\om\st \si\su x\}$$
for $\si\in 2^{<\om}$.  The following is the main lemma of
the proof of Theorem \ref{thm1}.

\begin{lemma} \label{main}
For each $s\in T$ not the root node and $\si\in 2^{<\om}$
$$(A_s\cap [\si]) \notin \mm_{<\al}$$
for $\al=\rank(s)$.
\end{lemma}
\proof
The proof is by induction on $\rank(s)$.  For $s\in \leaf(T)$,
i.e., $\rank(s)=0$, there is nothing to prove.
For $\rank(s)=1$ it easy to see by genericity that
$A_s$ is dense in $2^\om$ and so $A_s\cap[\si]$ cannot
be in $\mm_0$, the nowhere dense sets.

Working in $V$, for contradiction, choose $\al>1$ minimal so that
for some $s\in T$ with
$\rank(s)=\al$ there exists $p_0\in\poset$ and
$\si\in 2^{<\om}$ and $\be<\al$ such
that
$$p_0\forces (\name{A}_s\cap [\si])\in (\mm_\be)^\nn$$
Choose a hereditarily symmetric name $(\name{X}_n:n<\om)$ such that
$$p_0\forces\mbox{``} (\name{A}_s\cap [\si])
= \cup_{n<\om} \name{X}_n
\mbox{ where } \name{X}_n\in \mm_{\be_n}
\mbox{ for some } \be_n<\be<\al.\mbox{''}$$
Choose a finite $Q\su T$ such that
$H_Q$ fixes $\la \name{X}_n\st n<\om\ra$ and
$\dom(p_0)\su Q$.
Find an ordinal $\de$ with
\begin{enumerate}
 \item $\de\in\child(s)$,
 \item $\rank(\sde)\ge \be$, and
 \item $Q$ disjoint from $\{r\in T\st \sde\su r\}$.
\end{enumerate}
Choose an arbitrary $r\in \leaf(\sde)$.
Since $$ p_0\cup\{\la r,\si\ra\}\forces \name{x}_r\in \name{A}_s\cap [\si]$$
we can find an extension
$p_1\leq p_0\cup\{\la r,\si\ra\}$  and an $n_0$ so that
$$p_1\forces \name{x}_r\in \name{X}_{n_0}\cap[\si]$$
By extending $p_1$ even more, if necessary, we may assume
that $p_1(r)=\tau \supseteq \si$ where $\tau\in 2^{<\om}$
has the property that it is incompatible with $p_1(r^\pr)$
for every $r^\pr\in\dom(p_1)$ different from $r$.

\bigskip

Claim. $p_1\forces ([\tau]\cap \name{A}_{\sde})\su \name{X}_{n_0}$.

\bigskip\noindent
Suppose not. Then there exists $p_2\leq p_1$ and $r^\pr\supseteq \sde$
in $\dom(p_2)$ with $p_2(r^\pr)\supseteq \tau$ and
$$p_2\forces \name{x}_{r^\pr}\notin \name{X}_{n_0}.$$
Let $\pi\in \group$ be determined by the automorphism
of $\leaf(T)$ which swaps $r^\pr$ and $r$.
Note that $r^\pr\notin \dom(p_1)$ since $\tau$ was incompatible
with the range of $p_1$ except $p_1(r)$.  It follows from this that
$\pi(p_2)\cup p_1$ is a condition in $\poset$ (in fact
$\pi(p_2)\leq p_1$). By a general property
of automorphisms and forcing we have that
$$\pi(p_2)\forces \pi(\name{x}_{r^\pr})\notin \pi(\name{X}_{n_0}).$$
Since $\pi\in H_Q$ we have that $\pi(\name{X}_{n_0})=\name{X}_{n_0}$
and since $\hat{\pi}$ swaps $r^\pr$ and $r$ we have
that $\pi(\name{x}_{r^\pr})=\name{x}_r$ and so
$$\pi(p_2)\forces \name{x}_{r}\notin \name{X}_{n_0}.$$
But
$$p_1\forces \name{x}_{r}\in \name{X}_{n_0}$$
which contradicts the fact that $\pi(p_2)$ and $p_1$ are compatible.

The Claim contradicts the minimal choice of $\al$ since $\be_{n_0}<\al$
and $\mm_{\be_{n_0}}$ is closed under taking subsets.
This proves the lemma.
\qed

Working in $\nn$ for any ordinal $\al$ define
$\bb_\al$ to be all subsets of $2^\om$ whose symmetric difference
with an open set is in $\mm_\al$, i.e.,
$$\bb_\al=\{X\su 2^\om\st \exists U\su 2^\om \mbox{ open such that }
X\Delta U\in\mm_\al\}.$$

\begin{lemma} In the model $\nn$
$$\bsi(\al)\cup\bpi(\al)\su \bb_\al$$ for each $\al<\om_2$.
\end{lemma}
\proof
First we note that

\medskip\noindent
(a) $\bb_\al$ is closed under complementation.

\medskip
If $X\in\bb_\al$, then $(2^\om\sm X)\in \bb_\al$.
This is because, if $X=U\Delta Y$ where $U$ is open and
$Y\in\mm_\al$, then letting $Y^\pr=\cl(U)\sm U$, then
$Y^\pr\in\mm_0$ and so putting $V=2^\om\sm \cl(U)$ we have
that
$$(2^\om\sm X)\Delta V\su Y^\pr\cup Y\in \mm_\al.$$

Next we claim that

\medskip\noindent
(b)
If $\la X_n:n<\om\ra\in (\bb_{<\al})^\om$,
then $\cup_{n<\om}X_n\in\bb_{\al}$.

\medskip
We need to see we can
get the sequence of open sets required without using the axiom
of choice.

It follows from Lemma
\ref{main} that no nonempty open set is in $\mm_\al$ for
$\al<\om_2$.  An open set $U\su 2^\om$ is regular iff it is equal to the
interior of its closure, i.e., $U=int(cl(U))$.  If $U\su 2^\om$ 
is an arbitrary
open set, then $V=int(cl(U))$ is a regular open set containing $U$ such that
$V\Delta U$ is nowhere dense and hence in $\mm_0$.
($V\Delta U =V\sm U\su cl(U)\sm U$)

It follows that
for every $X\in\bb_\al$ there exists a regular open set $U$ such that
$X\Delta U\in\mm_\al$.

Suppose  $U$ and $V$ are regular open sets with $X\Delta U=A$ and
$X\Delta V=B$ where $A,B\in\mm_\al$.  Then
$U\Delta V = A\Delta B\su A\cup B\in\mm_\al$.   Since $\mm_\al$ contains
no nontrivial open sets and $U$ and $V$ are regular, it must be that $U=V$.

Hence for any $X\in \bb_\al$ there is a unique regular open
set $U$ such that $X\Delta U\in \mm_\al$.
Hence given $\la X_n:n<\om\ra\in (\bb_{<\al})^\om$,
choose $U_n$ the unique regular open set such that
$X_n\Delta U_n=Y_n\in \mm_{<\al}$.  Then
$$(\cup_{n<\om} X_n)\Delta(\cup_{n<\om}U_n)\su \cup_{n<\om}Y_n\in\mm_\al$$

From (a) and (b), induction and DeMorgan's Laws we have
that $\bpi(\al)$ and $\bsi(\al)$ are subsets of $\bb_\al$.

\qed

Next we prove the main Theorem of this paper.

\begin{theorem} \label{thm1}
It is relatively consistent with ZF that the Borel hierarchy on
$2^\om$ has length $\om_2$, i.e., the least $\al$ such that
$\bsi(<\al)$ is the
family of all Borel sets is $\al=\om_2$.
\end{theorem}
\proof
We show this holds in our model $\nn$.
Note
that if $\rank(s)=\al$ then $A_s\notin\bb_{<\al}$.  If it were,
then $A_s=U\Delta Y$ where $U$ open and $Y\in\mm_{<\al}$.
If $U$ is the empty set, then this would contradict Lemma \ref{main}.
But if $U$ is a nonempty set then $U\su A_s\cup Y$ and by
Lemma \ref{easy} $A_s\in\aa_\al\su\mm_\al$.  But Lemma
\ref{main} implies that no nontrivial open set is in $\mm_\al$.

It follows since each $A_s$ is Borel that the Borel hierarchy has
length at least $\om_2$.  But since $\om_2$ is a regular cardinal
in $\nn$ it must have length exactly $\om_2$.

\qed

Note that in $\nn$ if $X$ is any topological space which contains
a homeomorphic copy of $2^\om$, then the Borel order of $X$ is
$\om_2$.

\bigskip

Komj\'ath asks if 
the Borel hierarchy can have length greater than $\om_2$.
This would require a model in which
both $\om_1$ and $\om_2$ have cofinality $\om$.   In Gitik 1980
\cite{gitik} a model of ZF is produced (assuming the consistency
of ZFC plus unboundedly many
strongly compact cardinals) in which every $\aleph$
has cofinality $\om$.

In fact, we can prove

\begin{theorem} \label{gitikmodel}
Suppose $V$ is a countable transitive model of ZF in which
every $\aleph$ has countable cofinality.  Then for every
ordinal $\lambda$ in $V$, 
there is  symmetric submodel $\nn$
of a generic extension of $V$ with the same $\aleph$'s
as $V$ and the length of the
Borel hierarchy in $\nn$ is greater than $\lambda$.
\end{theorem}
\proof
We give a sketch of the proof at the end of this paper. 
\qed

\bigskip

\begin{center}
Countable unions of countable unions of etc., etc.
\end{center}

Specker 1957 \cite{specker} following Church 1927 \cite{church}
defines the classes  $\gg_\al$ for $\al$ an ordinal as follows:
\begin{enumerate}
\item $\gg_0$ is the class of countable sets,
\item $\gg_{<\al}=\cup_{\be<\al}\gg_\be$
\item $\gg_\al=\{\cup_{n<\om} X_n\st (X_n:n<\om)\in (\gg_{<\al})^\om\}$
\end{enumerate}
(Actually he defines $\gg_\al\sm\gg_{<\al}$.)
Gitik proves that in his model every set is in $\gg_{<\infty}$, i.e.,
$V=\gg_{<\infty}$.  L\"owe \cite{lowe} calls ZF+$V=\gg_{<\infty}$ the
theory ZFG and discusses some of its philosophical properties.

\begin{prop}
(Specker \cite{specker})
\begin{enumerate}
\item $\om_2$ is not the countable union of
countable sets, and in fact more generally
\item $\aleph_\al\notin \gg_{<\al}$ for any ordinal $\al$.  Similarly
\item $\pow(\aleph_\al)\notin \gg_\al$, and
\item if every $\aleph$ has cofinality $\om$, then
$\aleph_\al\in\gg_\al$ for every ordinal $\al$.
\end{enumerate}
\end{prop}
\proof
(1)  Suppose for contradiction that $\om_2=\cup_{n<\om}X_n$ where
each $X_n$ is countable.  For each $n<\om$ there exists a
unique countable ordinal $\al_n<\om_1$ and unique order preserving
bijection $f_n:\al_n\to X_n$.  Therefor there is no choice required
to define the onto map $f:\om\times\om_1\to \om_2$
by
$$f(n,\al)=\left\{
\begin{array}{ll}
f_n(\al) & \mbox{ if } \al<\al_n \\
0        & \mbox{ otherwise }
\end{array}\right.$$
But there is a definable bijection between $\om\times\om_1$
and $\om_1$ so this would be a contradiction.

(2) Left to the reader.

(3) In ZF there is a bijection between $\ka$ and $\ka\times\ka$ for any
infinite ordinal $\ka$.  Also there is a map from $\pow(\ka\times\ka)$
onto $\ka^+$ (map each well-ordering onto its order type).
Since $\gg_\al$ is closed under taking images and
$\aleph_{\al+1}\notin\gg_\al$ the result follows.

(4) $\aleph_0\in\gg_0$.  
Given $\aleph_\al$ we have by induction that
for every ordinal $\be<\aleph_\al$ that $\be\in\gg_{<\al}$ and
since the cofinality of $\aleph_\al$ is $\om$ the result follows.

\qed

It follows that in Gitik's model, $\om_2$ is the countable
union of countable unions of countable sets but cannot be
the countable union of countable sets.
In Gitik's model there is a simple example of a $\si$-algebra with a
long hierarchy:

\begin{prop}
Suppose every $\al\leq\om_2$ that $\cof(\aleph_\al)=\om$.
Let $\cc_0$ be the countable or co-countable subsets of $\aleph_{\om_2}$.
If $\cc$ is the $\si$-algebra generated by $\cc_0$,
then $\cc=\pow(\aleph_{\om_2})$ and
it takes exactly $\om_2+1$ steps to
generate $\cc$ from $\cc_0$ using
countable unions and countable intersections.
\end{prop}
\proof
$\aleph_{\om_2}\in\gg_{\om_2}\su
\cc$.  Since the $\gg$'s are closed under taking subsets,
We have that every subset of $\aleph_{\om_2}$ is in $\cc$.

Let $\comp(X)=\aleph_{\om_2}\sm X$ be the complement of $X$.
Define
$$\cc_\al=\{X\su \aleph_{\om_2}\st |X|\leq \aleph_\al \rmor
|\comp(X)|\leq\aleph_\al \}.$$
As usual $\cc_{<\al}=\cup_{\be<\al} C_\be$.
The following are easy to show:
\begin{enumerate}
\item $X\in\cc_\al$ iff $\comp(X)\in\cc_\al$.
\item If $\la X_n:n<\om\ra\ra\in (\cc_{<\al})^\om$, then
$\cup_{n<\om}X_n\in\cc_\al$ and $\cap_{n<\om}X_n\in\cc_\al$.
\item If $X\in\cc_\al$, then there exists
$\la X_n:n<\om\ra\ra\in (\cc_{<\al})^\om$ such either
$X=\cup_{n<\om}X_n$ or $X=\cap_{n<\om}X_n$.
\item If $A \su \aleph_{\om_2}$ has the property that
$|A|=|\comp(A)|=\aleph_{\om_2}$, then $A\notin \cc_{<\om_2}$.
\end{enumerate}
This shows that the hierarchy has exactly $\om_2+1$ levels.

\qed

A similar result holds for the sigma-field generated by
the countable subsets of $\aleph_{\om_3}$, etc.  Details are left to the
reader.

Unlike the $\aleph_\al$ the least $\ga$ such that
$\pow(\aleph_\al)$ gets into $\gg_\ga$ (if any) is not determined
by $\al$.  In the Feferman-Levy model
$\pow(\om)\in\gg_1\sm\gg_0$.  Gitik shows that in his model that
$\pow(\om)\in\gg_2\sm\gg_1$.  There is a variation of the Feferman-Levy
model where it is also true that $\pow(\om)\in\gg_2\sm\gg_1$.

We show that the least $\al$ such that
$\pow(\om)\in\gg_\al$ can be any $\al$ with $1\leq\al<\om_2$.
As in the proof of Theorem \ref{thm1} let $V$ be the Feferman-Levy model
and $T\in L$ be the well-founded tree of rank $(\aleph_{\om+1})^L$.
For each $\al<\om_2^V$ define
$$T_\al=\{s\st \alcats \in T\}.$$
Then the rank of $\la\ra$ in $T_\al$ is exactly
the rank of $\la\al\ra$ in $T$ which was $\al$.  Let
$\nn_\al$ be defined exactly as $\nn$ but using the tree
$T_\al$ in place of $T$.    Recall the definition of $\aa_\al$,
$\aa_0$ is the finite subsets of $2^\om$ and the $\aa_\al$ are
defined inductively as the countable unions of sets from
$\aa_{<\al}$.  So $\aa_{1+\al}$ is the same as $\gg_\al$ restricted
to subsets of $2^\om$.

\begin{theorem} \label{thm2}
For $2\leq\al<\om_2^V$
in the model $\nn_\al$, $\pow(\om)\in(\aa_\al\sm \aa_{<\al})$.
It follows that $\pow(2^\om)\su \aa_\al=$ Borel.  If $\al$ is a limit ordinal
then the Borel hierarchy in $\nn_\al$ has length exactly $\al$.
\end{theorem}

Only the statement $$(2^\om\in\aa_\al)^{\nn_\al}$$ needs to be proved. The
other parts of the Theorem are the same as Theorem \ref{thm1}. For example,
$\pow(2^\om)\su\aa_\al$, because  the $\aa_\al$ families are closed under taking
subsets.  $2^\om\notin \aa_{<\al}$ because the set
$A_{\la\ra}\notin\aa_{<\al}$.  Note that $\aa_{<\al}\su\mm_{<\al}$ and the rank
of $\la\ra$ in $T_\al$ is $\al$ (see Lemma \ref{main}).
The elements of $\aa_\al$ are
Borel because we started with finite sets and closed under taking countable
unions hence Borel $=\pow(2^\om)$.
If $\al$ is a limit ordinal then
$$\aa_{<\al}\su \bsi(<\al)\cap \bpi(<\al)\su\bb_{<\al}$$
Since the set
$A_{\la\ra}\notin\bb_{<\al}$, the Borel hierarchy has length
exactly $\al$.

\bigskip

The remainder of the proof of Theorem \ref{thm2}
(Lemmas \ref{clash}-\ref{compose}), is to
show that $2^\om\in\aa_\al$ holds in the model $\nn_\al$.
The intuitive reason this is true is because $A_{\la\ra}\in\aa_\al$ and
the reals in $\nn_\al$ can somehow be easily obtained from
$A_{\la\ra}$ and the reals in $V$.

Let $\la\cdot ,\cdot\ra$ be a recursive pairing function
from $\om\times\om$ to $\om$.  For example,
   $$\la n,m\ra=2^n(2m+1)-1$$
works.  Using this define a bijection
from $2^\om$ to $(2^\om)^\om$  by
$$x \mapsto (x_n\in 2^\om : n<\om) \mbox{ where } x_n(m)= x(\la n,m \ra).$$
Hopefully, we will not confuse the notation $x_n$ with the Cohen reals
$x_s$ which are attached to the nodes $s\in\leaf(T_\al)$.

For sets $A,B\su 2^\om$ define
$$A\clash B= \{x\in 2^\om\st
\exists N<\om\;\exists y\in B\;\;\forall n<N\; x_n\in A
\rmand
\forall n\geq N\;\; x_n=y_n\}$$

\begin{lemma} \label{clash}
For any $\al\geq 1$ if $A,B\in\aa_\al$, then
$A\clash B\in\aa_\al$.
\end{lemma}
\proof
For $\al=1$ note that for $A$ and $B$ countable, the
set $A\clash B$ is countable (without using choice).
Recall that the $\aa_\al$ families are closed under finite
unions.  Given increasing sequences $A_n$ and $B_n$ for
$n<\om$ note that
$$(\cup_{n<\om}A_n)\clash (\cup_{n<\om}B_n)=
\cup_{n<\om}(A_n\clash B_n)$$
So now the result follows by induction.
\qed

For $A\su 2^\om$ define
$$A^{<\om}=\{x\in 2^\om\st \exists N<\om\;\;
\forall n<N\; x_n\in A \rmand
\forall n\geq N\;\; x_n\equiv 0\}$$
where $x \equiv 0$ means $x$ is identically zero.

\begin{lemma}
For any $\al\geq 1$ if $A\in\aa_\al$, then
$A^{<\om}\in\aa_\al$.
\end{lemma}
\proof
 Note that $A^{<\om}=A\clash \{\underline{0}\}$ where
 $\underline{0}$ is the identically zero function.
\qed

In the model $V[G_{\al}]$ for each $t\in T_\al\sm\leaf(T_\al)$, define
$$B_t=\{x\in 2^\om: \exists s\supseteq t\;\;\rank(s)=1\;\;
\rmand \forall n<\om \;\;\;x_n=x_{\scatn}\}.$$
Recall that $A_t=\{x_s:s\in\leaf(t)\}$.
Define $C_t=A_t\clash B_t$.

\begin{lemma}
$C_t\in\nn_\al$, in fact, $C_t\in(\aa_{\be})^{\nn_\al}$ where
$\be=\rank(t)$.

\end{lemma}
\proof
Working in $V$ consider the set $P_t$ of sequences of
names, $\la \name{x}_n : n<\om\ra$ such that
there exists $N<\om$ and $s\supseteq t$ with
$\rank(s)=1$ such that
\begin{enumerate}
\item for all $n<N$ there exists  $r\in\leaf(t)$ such that
 $\name{x}_n =\name{x}_r$ and
\item for all  $n\geq N\;\;\name{x}_n=\name{x}_{\scatn}$.
\end{enumerate}
Recall that all $\pi\in\group$ have finite support and
the $\pi\in H_t$ permute the set of names for elements of $A_t$, i.e.,
$\{\name{x}_s:s\in\leaf(t)\}$,
moving only finitely many of them.  It follows that any
$\pi\in H_t$ permutes around
the elements of $P_t$.  From $P_t$ it is an exercise
to construct a name for $\name{C}_t$ which is fixed by
$H_t$.

But $\pi\in H_t$ also map $\name{A}_{\tde}$ to itself for
each $\de\in\child(t)$. Hence $H_t$ fixes the sequence
$(\name{C}_{\tde}:\de\in\child(t))$.
Recall that $\child(t)$ is countable in $V\su\nn_\al$ and
since
$$C_t=\bigcup\{ (\cup_{s\in F}A_s)\clash C_{\tde}\st
\de\in\child(t)\rmand F \in [\child(t)]^{<\om} \}$$
the lemma follows by induction.
\qed
Now we have by the Lemmas that since
$C_{\la\ra}\in\aa_\al$
in $\nn_\al$

\begin{cor} \label{cor}
$C^{<\om}_{\la\ra}\in\aa_\al$.
\end{cor}

Working in $V$ define $\qq$ to be the set of all
$f:\om\times\om\to 2^{<\om}\cup\{*\}$.
Since $\qq$ is
essentially the same as $\om^\om$ we know that
$\qq$ is the countable union of countable sets.
Given any $f\in\qq$ and $x\in 2^\om$
define $f(x)\in 2^\om$ by
$$f(x)(n)=\left\{
\begin{array}{ll}
1 & \mbox{ if }\exists m \;\; f(n,m)\su x \\
0 & \mbox{ otherwise }\\
\end{array}
\right.$$
We assume that $*$ is not a subsequence of any $x$.   For example,
if $M$ is a model of ZF and $x$ is $2^{<\om}$-generic over $M$, then
for any $y\in M[x]\cap 2^\om$ there exists $f\in M$ such that
$f(x)=y$.  To see this, work in $M$, and construct $f$ so
that  for any $n<\om$
$$\{f(n,m)\st m<\om\}=\{p\in 2^{<\om} \st
p\forces \name{y}(n)=1\}.$$

\begin{lemma}\label{char}
In $V[G]$, for all $y\in 2^\om$
$$y\in\nn_\al \rmiff \exists f\in \qq^V\;\exists z\in C^{<\om}_{\la\ra}
\;\; f(z)=y$$
\end{lemma}
\proof
The implication $\leftarrow$  is trivial because both
$\qq^V$ and  $C^{<\om}_{\la\ra}$ are in $\nn_\al$.

For the nontrivial direction,
we will find $z\in B^{<\om}_{\la\ra}$.
Suppose that
$y\in 2^\om\cap\nn_\al$ and suppose
$H_Q$ fixes $\name{y}$ where $Q$ is a finite subset of
$T_\al$.

At this point it would simplify our argument to assume that
for any $s\in T$ if $\rank(s)>1$ , then
the $\rank(\sde)>0$ for all $\de\in\child(s)$.
Equivalent, the parent of any leaf node has rank one.
Obviously we could have built $T$ with this property,
so we assume we did.

Assume that $Q$ contains the rank one parent
of every rank zero node in $Q$.
Let $(s_i:i<N)$ list all rank one nodes in $Q$.
Define 
\begin{enumerate}
\item $\leaf(Q)=\cup\{\leaf(s_i):i<N\}$ and
\item $\poset_Q=\{p\in\poset\st \dom(p)\su\leaf(Q)\}$.
\end{enumerate}
We claim that $y$ has a $\poset_Q$-name.  To see
this note that for any pair of finite sets $F_0$ and $F_1$ of leaf nodes
disjoint from $\leaf(Q)$ there is a $\pi\in H_Q$ for
which $\hat{\pi}(F_0)$ is disjoint from $F_1$.  From this it follows
that for any $n,i,$ and $p\in\poset$
$$p\forces \name{y}(n)=i \rmiff p\res_{\leaf(Q)}\forces \name{y}(n)=i.$$
Hence $y$ has a $\poset_Q$-name.

Define $z^i\in 2^\om$ for each $i<N$ so that $z^i_n=x_{\sin}$ for every $n$.
So
\begin{enumerate}
\item each $z^i$ is in $B_{\la\ra}$,
\item $y\in V[\la z_i:i<N\ra]$ and
\item $\la z_i:i<N\ra$ is $(2^{<\om})^N$-generic over $V$.
\end{enumerate}
As in the argument of Lemma \ref{regular}, let
$$A=\{(p,n,i)\in (2^{<\om})^N\times \om\times \{0,1\}
\st p\forces
\name{y}(n)=i\}.$$
Since there exists $n<\om$ with $A\in L[G_n]$, we can construct
$f\in L[G_n]\su V$ such that $f(\la z_i:i<N\ra)=y$.
\qed

\begin{lemma}\label{compose}
In $\nn$, for any set $A\in\aa_\al$ where $\al\geq 2$ the
set $$\qq\circ A = ^{def}\{f(x):f\in\qq\rmand x\in A\}$$ is in $\aa_\al$.
\end{lemma}
\proof
For $\al=2$ $\;\;\aa_\al$ is the family of sets which are the countable
union of countable sets. Let $A=\cup_n A_n$ and
let $\qq=\cup_n\qq_n$ where $A_n$ and $Q_n$ are countable.
Then for each $n,m<\om$ the set
$$\{f(x):x\in A_n\rmand f\in\qq_m\}$$
is countable, so $\qq\circ A$ is the countable union of
countable sets.

For larger $\al$ note that
$$\qq\circ (\cup_{n<\om}A_n)=\cup_{n<\om}\qq\circ A_n $$
so the result follows by induction.
\qed

\bigskip
\noindent
By Corollary \ref{cor} and Lemmas \ref{char} and \ref{compose},
we have that in $\nn_\al$
$$2^\om=\qq\circ C^{<\om}_{\la\ra}\in \aa_\al$$
hence this concludes the proof of Theorem \ref{thm2}.
\qed

\bigskip

{\bf Remark}. For successor ordinals $\al$ we get a
weaker result.
Suppose $\al=\lambda+n$ for $\lambda$ limit ordinal and $0<n<\om$, then
the Borel hierarchy in $\nn_\al$
has length $\ga$ where $\lambda+n\leq\ga\leq\lambda+2n$.
We are not sure what it is exactly.
The problem is that in the definition of $\bsi(\al)$
and $\bpi(\al)$
we forced an alternation between union and intersection.
Hence $$\aa_{\lambda+n}\su \bpi(\lambda+2n)\cap \bsi(\lambda+2n).$$
If instead we allow taking unions and then more unions, e.g.,
redefined
$\bsi(\al)$ (and similarly $\bpi(\al)$) as follows:
$$\bsi(\al)=\{\cup_{n<\om}A_n\st (A_n:n<\om)\in
(\bsi(<\al)\cup\bpi(<\al))^\om\}$$
then this problem disappears and the Borel hierarchy has length exactly
$\al$ even for successor ordinal case.

On the other hand, if we instead defined  ${\bsi(\al)}$ to be the
smallest class
of sets containing $\bpi(<\al)$ and closed under countable unions, then in
our models for Theorem \ref{thm2}, $\bsi(2)$ contains all subsets of $2^\om$.
Using a similar, alternative definition for $\bpi(\al)$, we can
get an alternative definition for the
length of the Borel hierarchy.

\begin{question}
Using this alternative definition of the length of the Borel hierarchy, can
it be greater than $\om_1$?
\end{question}

\begin{center}
The width of the Borel hierarchy
\end{center}

The Hausdorff terminology for the Borel hierarchy is
defined as follows:
 $F$ is the family of closed sets,
 $G$ is the family of open sets,
 $F_\si$ is the family of sets which can written as
the countable union of closed sets,
 $G_\de$ is the family of sets which can written as
the countable intersection of open sets,
 $F_{\si\de}$  is the family of sets which can written as
the countable intersection of $F_\si$ sets,
 etc.

In this terminology in the Feferman-Levy model
every subset of $2^\om$ is $F_{\si\si}$, since it
is the countable union of countable sets.  Hence
Borel $=F_{\si\si}=G_{\de\de}$.

\begin{prop}
(Without using the axiom of choice)
$F_{\si\de}\neq G_{\de\si}$ (equivalently
$\bpi(3)\neq\bsi(3)$).
\end{prop}
\proof
Let $\rat$ be the set of $x\in 2^\om$ which are eventually zero.
Define $P=\rat^\om\su (2^\om)^\om$.  We can identify $(2^\om)^\om$ with
$2^\om$ via a recursive pairing function as in the proof of
Theorem \ref{thm2}.  It is easy to check that $P$ is a
$F_{\si\de}$-set.  We show that $P$ cannot be $G_{\de\si}$.

\bigskip\noindent
{\bf Claim.} { \it Suppose $G\su (2^\om)^\om$ is a $G_\de$ set
and $(q_i\in\rat:i<n)$ has the property that
$$G\su \prod_{i<n}\{q_i\}\times \prod_{n\leq k<\om}\rat.$$
Then there exists $m>n$ and $(q_i\in\rat:n\leq i<m)$
such that
$$G\cap \left(\prod_{i<m}\{q_i\}\times \prod_{m\leq k<\om}\rat\right)
=\emptyset.$$
\bigskip

To prove the Claim assume for simplicity that $n=0$. So $G\su P$.  $G$ is not
dense else we could effectively construct $x\in G$ with the property that
$x_n\notin \rat$ for every $n$.  To see this write $G$ as a descending sequence
of dense open sets $U_n$ and construct sequences
$(s_n^m\in 2^{<\om}:m<N_n)$ with
\begin{enumerate}
\item $N_n<N_{n+1}<\om$,
\item $s^n_m\su s^{n+1}_m$ for $m<N_n$,
\item $\{x\in (2^\om)^\om\st \forall i<N_n\;\; s_i^n\su x_i\}\su U_n$, and
\item $s_m^{n+1}(k)=1$ for some $k>|s_m^n|$ and for all $m<N_n$.
\end{enumerate}
By taking the union of the $s_m^n$'s we get $x\in G$ such that
$x_n\notin \rat$ for all $n$.

Since $G$ is not dense it is easy to find the required $q_i$'s.
This proves the Claim.

\bigskip
Now we prove the Proposition.  Suppose for contradiction
$P=\cup_{n<\om}G_n$
where each $G_n$ is a $G_\de$.   Construct
$(q_i\in \rat:i<N_n)$ so that
$$G_n\cap \left(\prod_{i<N_n}\{q_i\}\times \prod_{N_n\leq k<\om}\rat\right)
=\emptyset$$
by applying the Claim to the $G_\de$ set
$$G_n\cap \left(\prod_{i<N_{n-1}}\{q_i\}\times
\prod_{N_{n-1}\leq k<\om}2^\om\right).$$
But then $(q_i:i<\om)\in P\sm \cup_{n<\om} G_n$ which is a contradiction.

\qed

\bigskip

Rather than using the terminology, $F_{\si\si\de\si\si}$, for example,
let us consider the following.
 For $f\in 2^{<\om_1}$ define the class $\Gamma_f$ as
follows:

\begin{enumerate}
\item $\Ga=\Ga_{\la\ra}$ be the family of clopen subsets
of $2^\om$
\item For $f:\de\to 2$ where $\de$ is a limit ordinal, define
$$\Ga_f=\cup\{\Ga_{f\res \al}\st \al<\de\}$$
\item For $f:\al+1\to 2$ define
$$\mbox{ if $f(\al)=0$ then }\Ga_f=
\{\cup_{n<\om}A_n\st (A_n:n<\om)\in
(\Ga_{f\res\al})^\om\}$$
$$\mbox{ if $f(\al)=1$ then }\Ga_f=
\{\cap_{n<\om}A_n\st (A_n:n<\om)\in
(\Ga_{f\res\al})^\om\}$$
\end{enumerate}

Hence $F_{\si\si\de\si\si}=\Ga_{\la 1,0,0,1,0,0\ra}$.

Note that $\Ga_{\la 0,0\ra}=\Ga_{\la 0\ra}=$ open sets and
$\Ga_{\la 1,1\ra}=\Ga_{\la 1\ra}=$ closed sets.  To rule
out these trivial collapses, we define nontrivial
$f:\de\to 2$ to be admissible if
$f(0)\neq f(1)$.

For $f$ and $g$ admissible define
$f\vleq g$ iff there exists a strictly increasing
$$\pi:\dom(f)\to\dom(g) \mbox { such that } \forall \al\in\dom(f)
\;\;f(\al)=g(\pi(\al)).$$   Note that if $f\vleq g$, then
$\Ga_f\su \Ga_g$.
Instead of looking for very long Borel hierarchies we can ask instead
for very wide Borel hierarchies:

\begin{conj}
It is relatively consistent with ZF that for every
$f$ and $g$ admissible
$$f\vleq g \rmiff \Ga_f\su \Ga_g.$$
\end{conj}

However, it is impossible that it be infinitely wide, by which we
mean:

\begin{prop}
For any infinite set $X$ of admissables there exists
distinct $f,g\in X$ with $f\vleq g$, hence $\Ga_f\su\Ga_g$.
\end{prop}
\proof

The ordering $\vleq$ is a well-quasiordering.  This is due to Nash-Williams
\cite{nash}.  We show how to avoid
using the axiom of choice.

A well-quasi ordering $(Q,\vleq)$ is a reflexive transitive relation such that
for every sequence $(f_n\st n<\om)\in Q^\om$  there exists
$n<m$ with $f_n\vleq f_m$.
Besides the fact that Nash-Williams proof may use the axiom of choice,
the set $X$ might be infinite but not contain an infinite sequence,
i.e., $X$ is Dedekind finite.

This particular quasi-ordering is absolute;  take $\pi$ witnessing
$f\vleq g$ by choosing the least possible value:
$$\pi(\al)=\min\be\geq\sup\{\pi(\ga)+1:\ga<\al\}\mbox{ such that }
f(\al)=g(\be).$$
If any $\pi$ works, the least possible value $\pi$ works.
It follows that for any two models $M\su N$ of set theory and
$f,g\in M$,
$$M\models f\vleq g \rmiff N\models f\vleq g$$
This is true even if $M$ and $N$ are nonwell-founded models.
To see that ZF proves our proposition, suppose not.  Then there
is a countable model $(M,E)$ of ZF which models
$M\models X$ is an infinite pairwise $\vleq$-incomparable family.
Using forcing we can generically add
a sequence $(f_n\in X\st n<\om^M)$ and get a model
$N\supseteq M$ which thinks there is an infinite sequence
($\om^N=\om^M$)
which is an $\vleq$-antichain.  But the inner model of $N$,
$((L[f_n\in X\st n<\om^N])^N,E^N)$, satisfies the axiom of choice
and hence the Nash-Williams Theorem is true, which is a contradiction.

\qed

\begin{center}
Arbitrarily long Borel hierarchies
\end{center}

We prove Theorem \ref{gitikmodel}.

Suppose $V$ is countable transitive model of ZF and $\lambda$ is an ordinal in
$V$. Suppose that in $V$ we have $\cof(\aleph_\ga)=\om$ for all $\ga<\lambda$. 
We find a symmetric submodel $\nn$ of a generic extension of $V$ with the same
$\aleph$'s as $V$ and the length of the Borel hierarchy in $\nn$ is at least
$\lambda$.

Let $\ka=\aleph_{\lambda}$ and
$$\poset=\{p:F\to 2^{<\om}\st F\in [\ka]^{<\om}\}.$$
For any $q=(X_n:n<\om)$ a partition of $\ka$ let
$$H_q=\{\pi\in\group\st \forall n\;\;\hat{\pi}(X_n)=X_n\}.$$
where $\group$ is the group of automorphisms of $\poset$ determined
by finite support permutations of $\ka$.  Take $\ff$ to
be the filter of subgroups determined by the set of all such $H_q$
and $\nn$ the symmetric model. 
Let $x_\al\in 2^\om$ be the Cohen real attached to
$\al$ and for $X\su\ka$ in $V$ let $A(X)=\{x_\al\st \al\in X\}$
in $V[G]$.

\begin{lemma} 
If $(X\in [\ka]^{\aleph_\al})^V$ and $\si\in 2^{<\om}$, then
$$\nn\models (A(X)\cap[\si])\notin \mm_{<\al}.$$
\end{lemma}
\proof
If $X$ is infinite, $A(X)$ is dense, so $A(X)\cap[\si]\notin\mm_0$ the 
nowhere dense sets.  

So suppose $\al>0$ and in $V$ write
$X$ as the disjoint union of sets $X_n$ for $n<\om$ of smaller
cardinality.  Suppose there exists $\be<\al$ and $p_0$ such that
$$p_0\forces A(X)\cap[\si]=\cup_n Y_n 
\mbox{ where } (Y_n:n<\om)\in(\mm_{<\be})^\om.$$
Suppose $H_q$ fixes the hereditarily symmetric names $(\name{Y}_n:n<\om)$.
By refining the $X_n$ and $q$ we may assume that
$q=(Z_n:n<\om)$ is a partition with $Z_{2n}=X_n$ for all $n$.
Choose $Z_{2n_0}$ with $|Z_{2n_0}|\geq\aleph_\be$
and disjoint
from the domain of $p_0$.  Choose an arbitrary $\de\in Z_{2n_0}$ and
find an extension $p_1\leq p_0\cup\{(\de,\si)\}$ and $n_1$ such that
$$p_1\forces x_\de\in Y_{n_1}.$$
Let $\tau=p_1(\al)$ and assume $\tau$ is incomparable with the other
elements of the range of $p_1$.

\bigskip

{\bf Claim.} $p_1\forces A(Z_{2n_0})\cap [\tau]\su Y_{n_1}$.

\medskip
Suppose not and take $p_2\leq p_1$ and
$\be\in Z_{2n_0}$ such that $p_2(\be)\supseteq \tau$ and
$$p_2\forces x_\be\notin Y_{n_1}.$$
Then the automorphism $\pi$ which swaps $\de$ and $\be$ is
in $H_q$ and fixes $\name{Y}_{n_1}$ but 
$p_1$ and $\pi(p_2)$ are compatible and
$\pi(p_2)\forces x_\de\notin Y_{n_1}$.  
\qed
The claim yields the Lemma.
\qed

Although we do not know if $V$ and $V[G]$ have the same cardinals,
we can show that $V$ and $\nn$ have the same cardinals.

\begin{lemma}
Suppose $f:\al\to\be$ be in $\nn$ where $\al$ and $\be$ are 
ordinal.  Then there exist in $V$ a countable $B\su\ka$ such
that $f\in V[G_B]$. 
\end{lemma}
\proof
Let $H_q$ fix $\name{f}$ where $q=(X_n:n<\om)$.   Let
$B=\cup\{X_n\st |X_n|<\om\}$.   Then $B$ is a countable subset
of $\ka$. By the usual automorphism argument 
$f\in V[G_B]$.
\qed

The partial order $\poset_B$ is countable in $V$ and so
$V$ and $V[G_B]$ have the same cardinals, i.e., if
$f:\ga\to\be$ is a map in $V[G_B]$, then in
$V$ there is map $g:\ga\times\om\to\be$ such that for
every $\de\;\;$
$f(\de)=g(\de,m)$ for some $m<\om$.

This finishes our sketch of the proof of Theorem \ref{gitikmodel}.  Note that to
use this method to get the Borel hierarchy to have length at least $\om_2+1$ 
requires $\om_2+1$ strongly compact cardinals.

\address


\newpage

\def\header{ Appendix \hfill Electronic version only\hfill }
\markboth\header\header

\setcounter{page}{1}
\setcounter{theorem}{0}

\begin{center}
 Appendix
\end{center}

\bigskip
The appendix is not intended for final publication but for
the electronic version only.

\bigskip

\begin{center}
Elementary forcing facts
\end{center}

Let $M$ be a countable transitive model of ZF.  Let $\poset$
be a partial order in $M$.  Define 
\begin{enumerate}
\item $G$ is a $\poset$-filter iff
   \begin{enumerate}
   \item $G\su\poset$
   \item $p\leq q$ and $p\in G$ implies $q\in G$
   \item $p,q\in G$ implies there exists $r\in G$ with $r\leq p$ and
   $r\leq q$.
   \end{enumerate}
\item $D\su\poset$ is dense iff for every $p\in\poset$ there
exists $q\leq p$ with $q\in D$.  
\item $G$ is $\poset$-generic over $M$ iff $G$ is
a $\poset$-filter and $G\cap D\neq\emptyset$ for every
$D\in M$ dense in $\poset$.

\item The $\poset$-names are defined inductively on rank.
$\tau$ is a $\poset$-name iff
each element of $\tau$ is of the form $(p,\si)$ where
$p\in\poset$ and $\si$ is a $\poset$-name.  

\item Given a $\poset$-filter
$G$ and $\poset$-name $\tau$, the realization of $\tau$ given
$G$ is defined inductively by
$$\tau^G=\{\si^G\st \exists p\in G\; (p,\si)\in\tau\}.$$

\item If $G$ is $\poset$-generic over $M$, then
$$M[G]=\{\tau^G\st \tau \mbox{ is a $\poset$-name in $M$}\}.$$

\item Forcing: $p\forces \theta(\vec{\tau})$ iff for every
$G$ $\poset$-generic over $M$ if $p\in G$ then
$M[G]\models\theta(\vec{\tau}^G)$.

\end{enumerate}

It is shown that if $M$ is a countable transitive model of $ZF$ then $M[G]$
is a countable transitive model of ZF with $M\su M[G]$.

This is proved using the two key properties of forcing:
\begin{enumerate}
\item (definability)
 For any formula $\theta(x_1,\ldots,x_n)$,
$$p\forces_{\poset}\theta(\tau_1,\ldots,\tau_n)$$
is definable in $M$ by a formula of the form
$\psi(p,\poset,\tau_1,\ldots,\tau_n)$.
\item (truth) If $M[G]\models\theta(\vec{\tau}^G)$, then
$$\exists p\in G\;\; p\forces \theta(\vec{\tau}).$$
\end{enumerate}

If $\pi$ is an automorphism of $\poset$ in $M$, then
$\pi$ extends to the $\poset$-names by induction on rank:
$$\pi(\tau)=\{(\pi(p),\pi(\si))\st (p,\si)\in\tau\}.$$

A basic fact about such automorphisms is

\begin{lemma} If $\pi$ is an automorphism of $\poset$ in $M$, then
for any formula $\theta$, $p\in\poset$, and
$\poset$-names, $\tau_1,\ldots,\tau_n$
$$p\forces\theta(\tau_1,\ldots,\tau_n) \rmiff
\pi(p)\forces\theta(\pi(\tau_1),\ldots,\pi(\tau_n)).$$
\end{lemma}
\proof
First prove by induction on rank that
$$\tau^{\pi^{-1}(G)}=\pi(\tau)^{G}$$
and note that $M[G]=M[\pi^{-1}(G)]$.

Then show that the following are equivalent:

\begin{enumerate}
\item $p\forces \theta({\tau})$.
\item For all $G$ $\poset$-generic over $M$ with $p\in G$
$M[G]\models \theta({\tau}^G)$.
\item For all $G$ $\poset$-generic over $M$ with $p\in \pi^{-1}(G)$
$M[\pi^{-1}(G)]\models \theta({\tau}^{\pi^{-1}(G)})$.
\item For all $G$ $\poset$-generic over $M$ with $\pi(p)\in G$
$M[G]\models \theta(\pi({\tau})^{G})$.
\item $\pi(p)\forces \theta({\pi(\tau)})$.
\end{enumerate}

We have written the parameters $\tau_1,\ldots,\tau_n$ as ${\tau}$
to shorten the notation.

\qed

\begin{center}
The symmetric submodel
\end{center}

Suppose that $\group$ is a group of automorphisms of $\poset$
in $M$.  Then we can define in $M$:

\begin{enumerate}
\item For any $\poset$-name $\tau$ the subgroup of $\group$:
$$\fix(\tau)=\{\pi\in\group\st \pi(\tau)=\tau\}.$$

\item  $\ff$ is a normal filter of subgroups of $\group$ iff
 \begin{enumerate}
 \item if $H\su K\su \group$ are subgroups and $H\in\ff$, then $K\in\ff$,
 \item if $H,K\in\ff$, then $H\cap K\in\ff$, and
 \item if $H\in\ff$ and $\pi\in\group$, then $\pi H\pi^{-1}\in \ff$.
 \end{enumerate}

\item $\tau$ is symmetric iff $\fix(\tau)\in\ff$.

\item  $\tau$ is hereditarily symmetric
iff $\tau$ is symmetric and $\si$ is
hereditarily symmetric for every $(p,\si)\in\tau$.

\end{enumerate}

Remark.
Suppose $H=\fix(\tau)$ and $\pi\in\group$. Then
$$\pi H\pi^{-1}\su \fix(\pi(\tau)).$$
Hence if $\tau$ is an hereditarily symmetric name and $\pi\in\group$
then $\pi(\tau)$ is an hereditarily symmetric name.

For $G$ which is $\poset$-generic over $M$ define the
symmetric model
$$\nn=\{\tau^G\st \tau \mbox{ is an hereditarily symmetric $\poset$-name 
in $M$ }\}.$$

\begin{theorem}\footnote{Jech \cite{jech}
assumes $M$ models AC. I don't know why.}
Suppose $M$ is a countable transitive model of
ZF.  In $M$, 
$\poset$ is a poset,
$\group$ is a subgroup of the automorphism
group of $\poset$, and $\ff$ is a normal filter.  Then
for any $G$ which $\poset$-generic over $M$, the symmetric
model $\nn$ is a transitive model of ZF such that
$M\su\nn\su M[G]$.
\end{theorem}
\proof
The fact that $\nn$ is transitive follows from the definition of
hereditarily symmetric names.
$M\su\nn$ because the canonical names
$$\check{x}=\{(1,\check{y})\st y\in x\}$$
are fixed by every automorphism of $\poset$.
$\nn\su M[G]$ is obvious.

\medskip\noindent
Axioms of ZF are true in $\nn$:
\begin{enumerate}

\item Pair.  A name for the pair
$\{\tau^G,\si^G\}$ is $\{(1,\tau), (1,\si)\}$ and
$$\fix(\tau)\cap \fix(\si)\su \fix(\{(1,\tau) (1,\si)\}.$$
It follows that if $\si$ and $\tau$ are hereditarily symmetric, then
so is this name for their pair.

\item Union.  Given $\name{x}$,  let
$$\name{y}=\{(p,\si)\st \exists (r,\rho)\in\name{x}\;\;
\exists s\;\; (s,\si)\in \rho\;\; p\leq s \wedge p\leq r\}$$
Then
$$\forces\name{y}=\cup\name{x}$$
and $\fix(\name{x})\su\fix(\name{y})$.
If $\name{x}$ is hereditarily symmetric, so is $\name{y}$.

\item Power Set. Given $\name{x}$ hereditarily symmetric,  let
$$Q=\{\si\st\exists p\in\poset\;\; (p,\si)\in \name{x}\}$$
each element of $Q$ is hereditarily symmetric.  Let
$$\name{y}=\{(p,\si)\st \si\su\poset\times Q\mbox{ is symmetric and }
p\forces \si\su\name{x}\}.$$
then $\name{y}$ is a hereditarily symmetric name for the power
set of $\name{x}$ in $\nn$. Note that the normality condition guarantees
that if $\si$ is hereditarily symmetric then so is $\pi(\si)$ for
every $\pi\in \group$.  Also if
$$p\forces \si\su\name{x}$$
and $\pi\in\fix(\name{x})$ then
$$\pi(p)\forces \pi(\si)\su\name{x}.$$
So $\fix(\name{x})\su\fix(\name{y})$.

\item Comprehension.  Given a formula $\theta(v,\vec{\tau})$ with
hereditarily symmetric parameters and a hereditarily symmetric
$\name{x}$ then defining $Q$ as before let
$$\name{y}=\{(p,\si)\in\poset\times Q\st
p\forces \si\in\name{x}\;\;\;\nn\models \theta(\si,\vec{\tau})\}.$$
If $\pi$ fixes $\name{x}$ and each $\tau_i$ then
$\pi(\name{y})=\name{y}$.

\item Replacement.  We may assume that $M$ is a definable
class in $M[G]$ by adding a predicate $\name{M}$
if necessary.  Since $M[G]$ models replacement and
$\nn$ is a definable class in $M[G]$ for any formula
$\theta(x,y)$ and set $A\in\nn$ there will be a
set $B\in M$ of hereditarily symmetric names
such that for every $a\in A$ if $\nn\models\exists y\;\theta(a,y)$
then there exist $\tau\in B$ such that $\nn\models \theta(a,\tau^G)$.
$$C=\{(1,\pi(\tau)):\tau\in B\rmand \pi\in\group\}$$
is hereditarily symmetric and $\{\tau^G:\tau\in B\}\su C^G\in\nn$.

\end{enumerate}

\qed

\begin{center}
The Feferman-Levy model
\end{center}

The Feferman-Levy Model $V$ is described in Jech \cite{jech}.  The
ground model satisfies $V=L$, let us call it $L$.
In $L$ let $\col$ be the following version of the 
Levy collapse of $\aleph_\om$:
$$\col=\{p:F\to \aleph_{\om}\st F\in[\om\times\om]^{<\om}\rmand
\forall (n,m)\in F\;\; p(n,m)\in\aleph_n\}.$$
The group
$\group$ of automorphisms of $\col$ are those which are
determined by finite support permutations of $\om\times\om$
which preserve the first coordinate,
that is, $\pi\in\group$ iff
there exists a finite support permutation
$\hat{\pi}:\om\times\om\to \om\times\om$
such that
$\hat{\pi}(n,m)=(n^\pr,m^\pr)$ implies $n=n^\pr$ and
$\pi(p)(s)=p(\hat{\pi}(s))$ for all $p\in\col$.
The normal filter $\ff$ of
subgroups is generated by
$$H_n=\{\pi\in\group\st
\hat{\pi}\res n\times\om \mbox{ is the identity
}\}$$
for $n<\om$.

The Feferman-Levy model, $V$, is the symmetric model
$L\su V\su L[G]$ determined
by $\col$, $G$, and the groups $\group, \ff$.

For any $n<\om$ let
$$\col_n=\{p\in\col\st \dom(p)\subseteq n\times\om\}.$$
For $G$ $\col$-generic over $L$ let $G_n=G\cap\col_n$.
Note that $H_n$ fixes the
canonical name for $G_n$,
$$\name{G}_n=\{(p,\check{p}): p\in\col_n\}$$
so $L[G_n]\su V$.  If we let
$$\name{X}_n=\{(1,\tau)\st \tau\su\col_n\times\{\check{k}:k<\om\}\}$$
then ${X}_n=L[G_n]\cap \pow(\om)$ and every $\pi\in\group$ fixes
$\name{X}_n$.  It follows that the sequence
$(L[G_n]\cap \pow(\om)\st n<\om)$ is in $V$.  Note that 
each $L[G_n]\cap \pow(\om)$ is countable in $V$.

\begin{theorem} \label{appendix3}
$$\pow(\om)\cap V=\bigcup_{n<\om}(L[G_n]\cap \pow(\om)).$$
More generally,
if $X\su Y\in L$ and $X\in V$, then for some $n<\om$ we have that
$X\in L[G_n]$
\end{theorem}

\proof
We prove the last statement.  Suppose
$$p_0\forces \name{X}\su\check{Y}\in L \rmand \name{X}\in V.$$
Choose $n$ large enough so that $H_n$ fixes $\name{X}$ and
$p_0\in\col_n$.

Note that for each $k\geq n$ that $\pi\in H_n$ can arbitrarily permute
$\{k\}\times\om$. It follows
that for any $y\in Y$ and $p\leq p_0$ that
$$p\forces \check{y}\in\name{X} \rmiff
   p\res_{(n\times\om)}\forces \check{y}\in\name{X}$$
and similarly
$$p\forces \check{y}\notin\name{X} \rmiff
   p\res_{(n\times\om)}\forces \check{y}\notin\name{X}.$$

Define
$$\name{W}=\{(p,\;\check{y})\in\col_n\times
\{\check{y}:y\in Y\}\st p\leq p_0
\rmand p\forces \check{y}\in \name{X}\}.$$
It follows that $p_0\forces \name{X}=\name{W}$.
But clearly, $W^G\in L[G_n]$.
\qed

\begin{center}
A variant of the Feferman-Levy model
\end{center}

We show that the following variant of the Feferman-Levy model
has the property that $\pow(\om)\in\gg_2\sm\gg_1$ using an argument
similar to Gitik's.   Redefine the Levy Collapse as follows:
$$\col=\{p:F\to \aleph_{\om}\st F\in[\aleph_\om\times\om]^{<\om}\rmand
\forall (\al,m)\in F\;\; p(\al,m)\in\al\}.$$
The group $\group$ is defined similarly, the normal filter
of subgroups, $\ff$, is defined to be the filter generated by
subgroups of the form
$$H_F=\{\pi\in\group\st \hat{\pi}\res F\times \om
\mbox{ is the identity} \}$$
where $F\in [\aleph_\om]^{<\om}$.
Call this alternative Feferman-Levy model $V^\pr$.

\begin{theorem}
In $V^\pr$ we have that $\pow(\om)$ is not the countable union
of countable sets but is the countable union of countable unions
of countable sets.
\end{theorem}
\proof

For any finite $F\su\aleph_\om$ define
$$\col_F=\{p\in \col\st \dom(p)\su F\times\om\}$$
and for $G$ which is $\col$-generic define
$$G_F=G\cap\col_F.$$

\bigskip
{\bf Claim.} $\pow(\om)\cap V^\pr=\cup\{L[G_F]\cap \pow(\om)\st F\in
[\om_1^V]^{<\om}\}.$

\bigskip This claim follows from a similar argument to the ordinary
Feferman-Levy model.

\medskip

Each $\col_F$-name is fixed by $H_F$.
The set of all $\col_F$-names:
$$\name{X}_F=\{(1,\tau): \tau \mbox{ is a $\col_F$-name}\}$$
is fixed by every $\pi\in\group$.
Note that $L[G_F]\cap \pow(\om)=X_F^G$ is a countable set in $V^\pr$
and the sequence $(X_F^G:F\in [\aleph_\om^L]^{<\om})$ is in $V^\pr$.
Note that
$$\bigcup_{n<\om}\cup\{L[G_F]\cap \pow(\om)\st F\in [\aleph_n^L]^{<\om}\}$$
is a countable union of countable unions of countable sets.

\medskip
Now we prove that in $V^\pr$ the
power set of $\om$
is not the countable union of countable sets.  This follows
from the

\bigskip
{\bf Claim.} If $Y\su X\in L$ and $Y\in V^\pr$,
then there exists $F$ finite such that $Y\in L[G_F]$.

This claim is proved similarly to Theorem \ref{appendix3}.

\medskip In $V^\pr$, suppose for contradiction
that $\pow(\om) =\cup_{n<\om} Y_n$ where each $Y_n$ is countable.
Working in $L$ let $(\name{Y}_n:n<\om)$ 
and $(\name{f}_n:n<\om)$ be sequences of hereditarily
symmetric names and $p\in\col$ such that for each $n$
$$p\forces \name{f}_n:\om\to\name{Y}_n\mbox{ is onto. }$$
By the Claim we can find in $L$ a sequence
$(F_n:n<\om)$ of finite sets such that
$$p\forces \name{f}_n\in L[G_{F_n}].$$ 
Choose any $\al\notin \cup_nF_n$ and let $x\su\om$ code
the generic map $g_\al:\om\to\al$.  Then $x\notin \cup_nY_n$.
\qed

\bigskip

\begin{center}
A remark on descriptive set theory
\end{center}

Levy \cite{levy} shows that in any model of ZF in which
$\om_1=\aleph_{\om}^L$ there is a $\Pi^1_2$ predicate
$Q(n,x)$ on $\om\times2^\om$ such that
$$\forall n\exists x\;Q(n,x)\;\; \wedge \;\; \neg\exists (x_n:n<\om)
\forall n \;Q(n,x_n).$$
The predicate $Q$ says that $x$ is a code for a countable model of the form
$(L_\al,\in)$ with $n$ infinite cardinals and there is no real $y$ coding a
model of the form $(L_\be,\in)$ with $\be>\al$ in which these cardinals are
collapsed.   He notes that such an example cannot be done for a
${\bf \Sigma}^1_2$ predicate because the Kondo-Addison Theorem can be proved
without the axiom of choice.

\bigskip

\begin{center}
Other interesting references.
\end{center}

Gregory H. Moore \cite{moore} has an interesting book on the history of the
axiom of choice.   H\'ajek \cite{hajek} shows the independence of Church's axioms
(although I have not been able to see a copy of this paper).  Hardy 1904
\cite{hardy,hardy2} shows  that $\om_1$ embeds into $\om^\om$ by building a strictly
increasing $\leq^*$ $\om_1$-sequence given a ladder sequence on $\om_1$, i.e.,
$(C_\al\su\al:\al^{lim}<\om_1)$ where $C_\al$ is a cofinal $\om$-sequence in
$\al$.


\begin{thebibliography}{99}

\bibitem{church} Church, Alonzo; Alternatives to Zermelo's
assumption.  Trans. Amer. Math. Soc.  29  (1927),  no. 1, 178--208.

\bibitem{gitik} Gitik, M.; All uncountable cardinals can be singular. Israel J.
Math. 35 (1980), no. 1-2, 61--88.

\bibitem{jech} Jech, Thomas J.; {\bf The axiom of choice.} Studies in Logic and
the Foundations of Mathematics, Vol. 75. North-Holland Publishing Co.,
Amsterdam-London; Amercan Elsevier Publishing Co., Inc., New York, 1973. xi+202
pp.

\bibitem{kechris} Kechris, Alexander S.; {\bf Classical descriptive set
theory.} Graduate Texts in Mathematics, 156. Springer-Verlag,
New York, 1995. xviii+402 pp.

\bibitem{leb} H.Lebesgue; Sur les fonctions repr\'{e}sentables  analytiquement,
Journal de Math\'{e}mati\-ques Pures et Appliqu\'{e}s, 1(1905), 139-216.

\bibitem{lowe} L\"owe, Benedikt;
A second glance at non-restrictiveness.
Philos. Math. (3) 11 (2003), no. 3, 323--331.

\bibitem{nash}
Nash-Williams, C. St. J. A.;
On better-quasi-ordering transfinite sequences.
Proc. Cambridge Philos. Soc. 64 1968 273--290.

\bibitem{specker} Specker, Ernst; Zur Axiomatik der Mengenlehre (Fundierungs-
und Auswahlaxiom). (German) Z. Math. Logik Grundlagen Math. 3 1957 173--210.

\end{thebibliography}

\begin{thebibliography}{99}

\bibitem{hajek} H\'ajek, Petr; The consistency of the Church's alternatives.
Bull. Acad. Polon. Sci. Sér. Sci. Math. Astronom. Phys. 14 1966 423--430.

\bibitem{hardy} Hardy, G.H.; A theorem concerning the infinite cardinal
numbers. Quarterly journal of pure and applied mathematics, 35(1904)
p.87-94.

\bibitem{hardy2} Hardy, G.H.; The continuum and the second number class,
Proceedings of the London mathematical society, (2) 4 (1906) p. 10-17.

\bibitem{levy} L\'evy, Azriel; Definability in axiomatic set theory. II. 1970
Mathematical Logic and Foundations of Set Theory (Proc. Internat. Colloq.,
Jerusalem, 1968) pp. 129--145 North-Holland, Amsterdam.

\bibitem{moore} Moore, Gregory H.; {\bf Zermelo's axiom of choice. Its origins,
development, and influence.} Studies in the History of Mathematics and Physical
Sciences, 8. Springer-Verlag, New York, 1982. xiv+410 pp. (1 plate). ISBN:
0-387-90670-3. 

\end{thebibliography}
\end{document}